\documentclass[twoside]{article}
\usepackage{tocstyle}

\usepackage[sc]{mathpazo}
\usepackage{graphicx,amsmath,amssymb,alltt}
\makeatletter
\providecommand*{\Dashv}{%
  \mathrel{%
    \mathpalette\@Dashv\vDash
  }%
}
\newcommand*{\@Dashv}[2]{%
  \reflectbox{$\m@th#1#2$}%
}
\makeatother

\usepackage[T1]{fontenc} 
\usepackage{microtype}
\usepackage{pgfplots,tikz,wrapfig}
\usetikzlibrary{cd}
\usepackage{lscape}

\usepackage{wasysym}
\usepackage{fancyvrb,color}

\newcommand{\underb}[1]{\tikz[baseline=(text.base)]{
    \node[inner sep=0pt,outer sep=0pt,text height = 2ex] (text) {#1};
    \draw[thin] (text.south east) -- ++(0,-0.1) -| (text.south west);}}

\def\@underlinebracket#1{\underline{\underb{#1}}}

\makeatletter
\def\PY@reset{\let\PY@it=\relax \let\PY@bf=\relax%
    \let\PY@ul=\relax \let\PY@tc=\relax%
    \let\PY@bc=\relax \let\PY@ff=\relax}
\def\PY@tok#1{\csname PY@tok@#1\endcsname}
\def\PY@toks#1+{\ifx\relax#1\empty\else%
    \PY@tok{#1}\expandafter\PY@toks\fi}
\def\PY@do#1{\PY@bc{\PY@tc{\PY@ul{%
    \PY@it{\PY@bf{\PY@ff{#1}}}}}}}
\def\PY#1#2{\PY@reset\PY@toks#1+\relax+\PY@do{#2}}

\expandafter\def\csname PY@tok@gd\endcsname{\def\PY@tc##1{\textcolor[rgb]{0.63,0.00,0.00}{##1}}}
\expandafter\def\csname PY@tok@gu\endcsname{\let\PY@bf=\textbf\def\PY@tc##1{\textcolor[rgb]{0.50,0.00,0.50}{##1}}}
\expandafter\def\csname PY@tok@gt\endcsname{\def\PY@tc##1{\textcolor[rgb]{0.00,0.27,0.87}{##1}}}
\expandafter\def\csname PY@tok@gs\endcsname{\let\PY@bf=\textbf}
\expandafter\def\csname PY@tok@gr\endcsname{\def\PY@tc##1{\textcolor[rgb]{1.00,0.00,0.00}{##1}}}
\expandafter\def\csname PY@tok@cm\endcsname{\let\PY@it=\textit\def\PY@tc##1{\textcolor[rgb]{0.25,0.50,0.50}{##1}}}
\expandafter\def\csname PY@tok@vg\endcsname{\def\PY@tc##1{\textcolor[rgb]{0.10,0.09,0.49}{##1}}}
\expandafter\def\csname PY@tok@vi\endcsname{\def\PY@tc##1{\textcolor[rgb]{0.10,0.09,0.49}{##1}}}
\expandafter\def\csname PY@tok@vm\endcsname{\def\PY@tc##1{\textcolor[rgb]{0.10,0.09,0.49}{##1}}}
\expandafter\def\csname PY@tok@mh\endcsname{\def\PY@tc##1{\textcolor[rgb]{0.40,0.40,0.40}{##1}}}
\expandafter\def\csname PY@tok@cs\endcsname{\let\PY@it=\textit\def\PY@tc##1{\textcolor[rgb]{0.25,0.50,0.50}{##1}}}
\expandafter\def\csname PY@tok@ge\endcsname{\let\PY@it=\textit}
\expandafter\def\csname PY@tok@vc\endcsname{\def\PY@tc##1{\textcolor[rgb]{0.10,0.09,0.49}{##1}}}
\expandafter\def\csname PY@tok@il\endcsname{\def\PY@tc##1{\textcolor[rgb]{0.40,0.40,0.40}{##1}}}
\expandafter\def\csname PY@tok@go\endcsname{\def\PY@tc##1{\textcolor[rgb]{0.53,0.53,0.53}{##1}}}
\expandafter\def\csname PY@tok@cp\endcsname{\def\PY@tc##1{\textcolor[rgb]{0.74,0.48,0.00}{##1}}}
\expandafter\def\csname PY@tok@gi\endcsname{\def\PY@tc##1{\textcolor[rgb]{0.00,0.63,0.00}{##1}}}
\expandafter\def\csname PY@tok@gh\endcsname{\let\PY@bf=\textbf\def\PY@tc##1{\textcolor[rgb]{0.00,0.00,0.50}{##1}}}
\expandafter\def\csname PY@tok@ni\endcsname{\let\PY@bf=\textbf\def\PY@tc##1{\textcolor[rgb]{0.60,0.60,0.60}{##1}}}
\expandafter\def\csname PY@tok@nl\endcsname{\def\PY@tc##1{\textcolor[rgb]{0.63,0.63,0.00}{##1}}}
\expandafter\def\csname PY@tok@nn\endcsname{\let\PY@bf=\textbf\def\PY@tc##1{\textcolor[rgb]{0.00,0.00,1.00}{##1}}}
\expandafter\def\csname PY@tok@no\endcsname{\def\PY@tc##1{\textcolor[rgb]{0.53,0.00,0.00}{##1}}}
\expandafter\def\csname PY@tok@na\endcsname{\def\PY@tc##1{\textcolor[rgb]{0.49,0.56,0.16}{##1}}}
\expandafter\def\csname PY@tok@nb\endcsname{\def\PY@tc##1{\textcolor[rgb]{0.00,0.50,0.00}{##1}}}
\expandafter\def\csname PY@tok@nc\endcsname{\let\PY@bf=\textbf\def\PY@tc##1{\textcolor[rgb]{0.00,0.00,1.00}{##1}}}
\expandafter\def\csname PY@tok@nd\endcsname{\def\PY@tc##1{\textcolor[rgb]{0.67,0.13,1.00}{##1}}}
\expandafter\def\csname PY@tok@ne\endcsname{\let\PY@bf=\textbf\def\PY@tc##1{\textcolor[rgb]{0.82,0.25,0.23}{##1}}}
\expandafter\def\csname PY@tok@nf\endcsname{\def\PY@tc##1{\textcolor[rgb]{0.00,0.00,1.00}{##1}}}
\expandafter\def\csname PY@tok@si\endcsname{\let\PY@bf=\textbf\def\PY@tc##1{\textcolor[rgb]{0.73,0.40,0.53}{##1}}}
\expandafter\def\csname PY@tok@s2\endcsname{\def\PY@tc##1{\textcolor[rgb]{0.73,0.13,0.13}{##1}}}
\expandafter\def\csname PY@tok@nt\endcsname{\let\PY@bf=\textbf\def\PY@tc##1{\textcolor[rgb]{0.00,0.50,0.00}{##1}}}
\expandafter\def\csname PY@tok@nv\endcsname{\def\PY@tc##1{\textcolor[rgb]{0.10,0.09,0.49}{##1}}}
\expandafter\def\csname PY@tok@s1\endcsname{\def\PY@tc##1{\textcolor[rgb]{0.73,0.13,0.13}{##1}}}
\expandafter\def\csname PY@tok@dl\endcsname{\def\PY@tc##1{\textcolor[rgb]{0.73,0.13,0.13}{##1}}}
\expandafter\def\csname PY@tok@ch\endcsname{\let\PY@it=\textit\def\PY@tc##1{\textcolor[rgb]{0.25,0.50,0.50}{##1}}}
\expandafter\def\csname PY@tok@m\endcsname{\def\PY@tc##1{\textcolor[rgb]{0.40,0.40,0.40}{##1}}}
\expandafter\def\csname PY@tok@gp\endcsname{\let\PY@bf=\textbf\def\PY@tc##1{\textcolor[rgb]{0.00,0.00,0.50}{##1}}}
\expandafter\def\csname PY@tok@sh\endcsname{\def\PY@tc##1{\textcolor[rgb]{0.73,0.13,0.13}{##1}}}
\expandafter\def\csname PY@tok@ow\endcsname{\let\PY@bf=\textbf\def\PY@tc##1{\textcolor[rgb]{0.67,0.13,1.00}{##1}}}
\expandafter\def\csname PY@tok@sx\endcsname{\def\PY@tc##1{\textcolor[rgb]{0.00,0.50,0.00}{##1}}}
\expandafter\def\csname PY@tok@bp\endcsname{\def\PY@tc##1{\textcolor[rgb]{0.00,0.50,0.00}{##1}}}
\expandafter\def\csname PY@tok@c1\endcsname{\let\PY@it=\textit\def\PY@tc##1{\textcolor[rgb]{0.25,0.50,0.50}{##1}}}
\expandafter\def\csname PY@tok@fm\endcsname{\def\PY@tc##1{\textcolor[rgb]{0.00,0.00,1.00}{##1}}}
\expandafter\def\csname PY@tok@o\endcsname{\def\PY@tc##1{\textcolor[rgb]{0.40,0.40,0.40}{##1}}}
\expandafter\def\csname PY@tok@kc\endcsname{\let\PY@bf=\textbf\def\PY@tc##1{\textcolor[rgb]{0.00,0.50,0.00}{##1}}}
\expandafter\def\csname PY@tok@c\endcsname{\let\PY@it=\textit\def\PY@tc##1{\textcolor[rgb]{0.25,0.50,0.50}{##1}}}
\expandafter\def\csname PY@tok@mf\endcsname{\def\PY@tc##1{\textcolor[rgb]{0.40,0.40,0.40}{##1}}}
\expandafter\def\csname PY@tok@err\endcsname{\def\PY@bc##1{\setlength{\fboxsep}{0pt}\fcolorbox[rgb]{1.00,0.00,0.00}{1,1,1}{\strut ##1}}}
\expandafter\def\csname PY@tok@mb\endcsname{\def\PY@tc##1{\textcolor[rgb]{0.40,0.40,0.40}{##1}}}
\expandafter\def\csname PY@tok@ss\endcsname{\def\PY@tc##1{\textcolor[rgb]{0.10,0.09,0.49}{##1}}}
\expandafter\def\csname PY@tok@sr\endcsname{\def\PY@tc##1{\textcolor[rgb]{0.73,0.40,0.53}{##1}}}
\expandafter\def\csname PY@tok@mo\endcsname{\def\PY@tc##1{\textcolor[rgb]{0.40,0.40,0.40}{##1}}}
\expandafter\def\csname PY@tok@kd\endcsname{\let\PY@bf=\textbf\def\PY@tc##1{\textcolor[rgb]{0.00,0.50,0.00}{##1}}}
\expandafter\def\csname PY@tok@mi\endcsname{\def\PY@tc##1{\textcolor[rgb]{0.40,0.40,0.40}{##1}}}
\expandafter\def\csname PY@tok@kn\endcsname{\let\PY@bf=\textbf\def\PY@tc##1{\textcolor[rgb]{0.00,0.50,0.00}{##1}}}
\expandafter\def\csname PY@tok@cpf\endcsname{\let\PY@it=\textit\def\PY@tc##1{\textcolor[rgb]{0.25,0.50,0.50}{##1}}}
\expandafter\def\csname PY@tok@kr\endcsname{\let\PY@bf=\textbf\def\PY@tc##1{\textcolor[rgb]{0.00,0.50,0.00}{##1}}}
\expandafter\def\csname PY@tok@s\endcsname{\def\PY@tc##1{\textcolor[rgb]{0.73,0.13,0.13}{##1}}}
\expandafter\def\csname PY@tok@kp\endcsname{\def\PY@tc##1{\textcolor[rgb]{0.00,0.50,0.00}{##1}}}
\expandafter\def\csname PY@tok@w\endcsname{\def\PY@tc##1{\textcolor[rgb]{0.73,0.73,0.73}{##1}}}
\expandafter\def\csname PY@tok@kt\endcsname{\def\PY@tc##1{\textcolor[rgb]{0.69,0.00,0.25}{##1}}}
\expandafter\def\csname PY@tok@sc\endcsname{\def\PY@tc##1{\textcolor[rgb]{0.73,0.13,0.13}{##1}}}
\expandafter\def\csname PY@tok@sb\endcsname{\def\PY@tc##1{\textcolor[rgb]{0.73,0.13,0.13}{##1}}}
\expandafter\def\csname PY@tok@sa\endcsname{\def\PY@tc##1{\textcolor[rgb]{0.73,0.13,0.13}{##1}}}
\expandafter\def\csname PY@tok@k\endcsname{\let\PY@bf=\textbf\def\PY@tc##1{\textcolor[rgb]{0.00,0.50,0.00}{##1}}}
\expandafter\def\csname PY@tok@se\endcsname{\let\PY@bf=\textbf\def\PY@tc##1{\textcolor[rgb]{0.73,0.40,0.13}{##1}}}
\expandafter\def\csname PY@tok@sd\endcsname{\let\PY@it=\textit\def\PY@tc##1{\textcolor[rgb]{0.73,0.13,0.13}{##1}}}

\makeatother

\usepackage[hmarginratio=1:1,top=32mm,columnsep=20pt]{geometry} 
\usepackage{multicol} 
\usepackage[hang, small,labelfont=bf,up,textfont=it,up]{caption} 
\usepackage{booktabs} 
\usepackage{float} 
\usepackage[hidelinks]{hyperref} 
\usepackage{lettrine} 
\usepackage{paralist} 
\usepackage{abstract} 
\usepackage{outlines} 
\usepackage{mathtools}

\usepackage{titlesec} 
\renewcommand\thesection{\Roman{section}} 
\renewcommand\thesubsection{\Roman{subsection}} 
\titleformat{\section}[block]{\large\scshape\centering}{\thesection.}{1em}{} 
\titleformat{\subsection}[block]{\large}{\thesubsection.}{1em}{} 

\usepackage{fancyhdr} 
\title{\vspace{-15mm}{\large{
\selectfont\textbf{The Constituents of Sets, Numbers, and Other Mathematical Objects\\
Part One}}}}

\small
\author{
\textsc{\small Ruadhan O'Flanagan}\footnote{\href{mailto:rof@ruadhan.net}{rof@ruadhan.net}}
\vspace{-3mm}
}
\vspace{3mm}
{\small\date{\small May 2019}}
\usepackage[format=plain,labelsep=newline,singlelinecheck=false]{caption}
\begin{document}

\maketitle 

\thispagestyle{fancy} 

\vspace{-5mm}

\begin{abstract}
\hrule
\vspace{0.5mm}
\hrule

\vspace{1mm}

\noindent

\small 

The sets used to construct other mathematical objects are pure sets, which means that all of their elements are sets, which are themselves pure. One set may therefore be within another, not as an element, but as an element of an element, or even deeper, inside several layers of sets within sets. 

The introduction of the term {\em constituent} to describe a set which is within a given set, however deep, induces an apparently novel partial order on sets, and assigns to any given set a diagram which specifies a directed graph, or category, herein dubbed its constituent structure, indicating which sets within it are constituents of which others. 

Sets with different numbers of elements can have exactly the same constituent structure. Consequently, constituent structure isomorphisms between sets need not preserve the number of elements, although they are still injective, surjective, and invertible. We consider in detail an example of an isomorphism between a one-element set and a five-element set, which is a surjective mapping despite the mismatch in cardinalities.

The constituent structure of a set determines the mathematical objects for which the set is a suitable representation. Different schemes for constructing the natural numbers, such as those of von Neumann and Zermelo, generate sets with the same constituent structures. Objects share the constituent structures, not the elements, of the sets used to construct or represent them.   

The requirement that an object's properties be faithfully encoded within a set's constituent structure and not its non-constituent characteristics such as its cardinality, when made explicit, dictates a specific and novel way of representing ordered pairs and tuples of sets as sets, providing simple formulae for addressing and extracting sets located deep within nested tuples.

\vspace{2mm}

\hrule
\vspace{0.5mm}
\hrule
\end{abstract}

\section*{Introduction}

The subset relation, $\subset$, has the properties that make life easy and pleasant for those who wish to study sets. It is transitive, $a\subset b {\ \& \ } b\subset c\implies a\subset c$, which provides plenty of opportunities to draw deductions from given information, and it is antisymmetric, $a\subset b {\ \& \ } b\subset a\implies a=b$, which makes it possible for us to figure out what a set that we are interested in is equal to, which is generally at least a minor success.

The membership relation, $\in$, however, is not so pleasant to work with. It is not transitive, the sets it relates are never equal, and given a list of sets, $A, B, C, \cdots$, and membership relations between them, $B\in C, A\in C, \cdots$, there is no way of combining the given statements that permits us to draw a conclusion.

A careful comparison of the two basic relations between sets reveals that we don't like membership very much and would like it to go away. The standard ridding procedure is to contend that membership doesn't express any information that isn't accessible through the subset relation, which makes membership a superfluous and irrelevant part of the formalism. We can define membership using the subset relation: $A\in B \iff \{A\}\subset B$, so subsets give us everything we need, we tell ourselves.

Unfortunately, when we do this, we lose the ability to reference the structure deep inside sets. 

When we construct mathematical objects of interest, it typically involves gathering certain previously defined objects, $x_1, x_2, \cdots$, into a set, $X$. Then we equip that set with some structure, $S$, such as a function or a relation, by putting $X$ and $S$ together into an ordered pair, $P=(X, S)$, which is a specific type of set. We might combine two objects of that type, $P_1$ and $P_2$, with another structure, such as an isomorphism, $M$, to get an object, $K$, that we want to talk about, $K=(P_1, P_2, M)$.

The result is a set, $K$, which, deep inside, has sets whose internal structures replicate the properties of the objects $x_1, x_2, \cdots$, and larger sets which specify additional structure.

If somebody else constructs that set, $K$, and then gives it to us, and we try to understand what it is, using all the information that can be expressed in terms of the subset relation, so that we can list, as distinct symbols, all of the subsets of $K$, and say which of those are subsets of which others, we will be able to deduce that there are three things inside $K$ and nothing further.

This disappointing achievement is formalized on a grand scale in category theory, in which sets are regarded as equivalent if there is an invertible mapping, or isomorphism, between them. These mappings are assumed to be functions, which assign to each element of one set a corresponding element of the other. A one-to-one correspondence of elements, which establishes that two sets have the same number of elements, renders the sets equivalent, and, within the category of sets, indistinguishable.

Despite this drastic loss of information about the contents of a set, category theory appears to be the right way to represent structure, effectively using directed graphs instead of sets as the most elementary objects, thereby putting relations on the same footing as the objects they relate, instead of building an object to encode a relation and placing that relation object into a set along with the objects it relates.

A category is a neatly self-contained chunk of structure which can be isolated from other structures and studied in isolation. The size, or cardinality, of finite sets provides an obvious example of information that can be studied on its own while everything else about those sets is forgotten. This is the information which is captured and probed by the category of sets.

This leads to the questions: What category contains the other information about the set, apart from its cardinality? What structure is left when a set's cardinality is taken from it?

We can access information deep inside a set using the membership relation, $\in$, but we can't use that relation to specify which dots in a graph depicting a category should have arrows connecting them, because a relation in a category must be reflexive, which means that every dot in the graph must have an arrow from itself to itself, and relations depicted by arrows must be transitive. $\in$ is neither reflexive nor transitive.

We also know that the $\in$ relation specifies all of the information inside a pure set. If we know its elements, and their elements and so on, all the way down to the empty set, then we know everything about the pure set. But part of this information, namely cardinality, is isolated by the category of sets, and when that's removed, the remainder must be less informative than $\in$.

So we seek another relation, which is transitive and reflexive, and discards some information.
That search leads to the constituency relation, and to the category of constituent structures.

\subsection*{Outline}

The sets studied in this paper and its sequel in Part Two are all finite pure sets. A theory of the constituents of infinite sets requires a separate investigation to be undertaken after the finite case has been fully understood with clarity and certainty. For this reason, the construction of numbers ends with 
the rationals, since the construction of the real numbers requires the use of infinity.

Because every set considered is both pure and finite, it can be constructed in a finite number of steps, starting from the empty set and successively gathering previously generated sets together in a new set. There is consequently no need or opportunity to worry about whether the sets exist or are compatible with any specific system of axioms, other than those of pure finite sets.

In the first section, immediately following this, the constituents of a set and 
its constituent structure, along with the corresponding diagram, are defined and described. 

Section 2 uses the different constructions of the natural numbers introduced by von Neumann and Zermelo as examples of constituent structures which are shared by apparently very different sets. The isomorphism between Zermelo's number 5 and von Neumann's number 5 is expounded as an example of an invertible mapping between sets with different numbers of elements, and a surjection from a one-element set to a five-element set.

Some operations on sets relevant to constituent structure are introduced and given notation in the section after that, and then section 4 shows how set operations relate to natural numbers in specific representations as sets. Subsequent reflections reveal the guiding principle for the construction of mathematical objects using sets, namely that the properties of objects should be encoded only within the constituent structure of the underlying sets.

The next section examines the representation within constituent structures of ordered pairs and tuples of sets. Formulae are developed for placing a set at a position within another, asserting that a set has something at a position or that one set is at a position within another, and extracting a set from a position arbitrarily deep inside tuples of tuples.

We then combine the ability to position arbitrary sets within a constituent structure with the operations on sets defined in section 3. This allows us to connect complex constituent structures together by fusing them at multiple vertices.

This completes Part One of this two-part series. In Part Two, we consider the natural representation within sets, and their constituent structure, of integers, arithmetic expressions, and rational numbers, which leads to new and unexpected insights.

\newpage
\tableofcontents
\newpage

\section{The Constituents of a Set}

\subsection{Definition and Notation}

The foundational concept is most easily understood when the definition is given in words:

\vspace{3mm}

{\bf Definition:} The constituents of a set are the set itself and the constituents of its elements.

\vspace{3mm}

To compile a list of a set's constituents, we first add the set itself to the list. Then we add the constituents of each of the set's elements. Those include, for each element, the element itself, and the constituents of its elements. So the set, its elements, their elements, and so on, are all included.

The choice of word is appropriate because this formal definition coincides with the existing usage of the verb, ``to constitute'': A thing constitutes itself, and, at the same time, its parts constitute it.

In symbols, the definition has the form:
\begin{align}
x\lhd y\  \equiv\  x=y \text{\sc\small\ or } x\lhd z \in y
\end{align}
where $\lhd$ is the symbol for ``is a constituent of'' and $\equiv$ indicates that the left hand side is defined by the right hand side. The condition $x\lhd z \in y$ means that {\em some} element, $z$, of $y$, has $x$ as a constituent.

The symbolic definition can be read in words as: ``'$x$ is a constituent of $y$' is defined to mean that $x$ is equal to $y$ or $x$ is a constituent of some element, $z$, of $y$''.

The symbol $\lhd$ was chosen so that the expression $x\lhd y$ can be interpreted to graphically depict $x$ as something small and pointlike at the leftmost tip of the triangle, which is included within the larger object $y$ depicted by the side of the triangle opposite to $x$.

\vspace{2mm} 

\subsection{Properties}

Constituency has the useful properties which the membership relation lacked.

\vspace{2mm} 
\noindent $\bullet$ Transitivity
\vspace{1mm} 
\nopagebreak

Each set in a membership chain such as $A\in B\in C\in D$ is an element of the next set, but is a constituent of every set farther down the chain. In this example, the sets satisfy $A\lhd B\lhd C\lhd D$ automatically, and from these relations we can deduce others such as $A\lhd C$ and so on. So we can draw conclusions from multiple instances of the $\in$ relation after all, but only if we have the $\lhd$ symbol to express them.

\vspace{2mm} 
\noindent $\bullet$ Reflexivity
\vspace{1mm} 
\nopagebreak

The fact that a set is a constituent of itself distinguishes constituency from a previous way of constructing a form of ``transitive membership''. That way involved the concept of a transitive set, which is a set that contains the elements of its elements. Obviously, the set itself could not be included among those.

The reflexive property, $x\lhd x$, is central to the concept of constituency. The constituency relation with reflexivity provides sets with the structure of a category, which essentially means that something has clicked into place. A large body of mathematical knowledge now automatically applies to collections of sets and their constituency relations.

{
\nopagebreak
\vspace{2mm} 
\noindent $\bullet$ Antisymmetry
\vspace{1mm} 

Our attention is restricted here to pure finite sets. Two such sets, $x$ and $y$, which satisfy $x\lhd y \lhd x$, must clearly be equal to one another, because $x$ cannot be inside a smaller part of itself.

However, it is worth briefly observing that without the restriction to pure finite sets, a proof that $x\lhd y \lhd x \iff x=y$ would require the axiom of regularity, which states that every non-empty set contains an element that is disjoint from itself. The truth of that axiom is not obvious even for finite sets, but the possibility of membership loops among sets can be eliminated with the more obvious statement that no set is a constituent of any of its elements, which is equivalent to the statement that constituency is antisymmetric.
}

\vspace{2mm} 
\noindent $\bullet$ Partial Order
\vspace{1mm} 
\nopagebreak

The three properties of transitivity, reflexivity and antisymmetry make constituency a partial order on sets. 
In that context, there is a least element - the empty set, $\varnothing\equiv\{\}$, which satisfies $\varnothing\lhd x$ for every set, $x$.

When we consider only the constituents of a specific set, $y$, there is also a greatest element, namely the set $y$ itself. The ordering of the set's constituents can be shown in a diagram, with $y$ at the top and the empty set at the bottom, revealing the internal structure of the set.

\subsection{Constituent Structure}

\label{constituentstructure}

Formally, the constituent structure of a set is the category whose objects are the set's constituents and whose morphisms are constituency relations.

Informally, the constituent structure of a set is the specification of which of its constituents are contained as constituents in which others. 

We can display this structure for a set, $S$, in a diagram consisting of distinct horizontal levels, constructed using the following procedure:

\begin{enumerate}
\item
    Assign unique symbols, $A, B, \cdots$, to each set which is a constituent of $S$. Sets which occur more than once get only one symbol - symbols are assigned to sets, not occurrences of sets.
\item
    Add the symbol $S$ to the diagram. This is the top horizontal level.
\item
    For each symbol, $X$, at the lowest horizontal level, add the symbols of that set's elements one level below $X$. 
\item
    For any symbol which occurs more than once in the diagram, remove all occurrences apart from a single instance at the lowest level.
\item
    Repeat steps 3 and 4 until every constituent has been added to the diagram exactly once. The bottom horizontal level will contain just one symbol indicating the empty set.
\item
    Draw one edge for each membership relation, connecting the symbol of a set to the symbol of each set which contains it as an element.
\item
    Remove any edges connecting symbols which are connected by an upward path of two or more edges. These are the edges showing membership, $A\in C$, for which a membership chain, $A\in B\in C$, exists.
\end{enumerate}

Technically, this is the Hasse diagram of the set's constituents ordered by constituency. 

\begin{wrapfigure}{R}{3cm}
\centering
\begin{tikzcd}
 S \ar[-]{d}\\
B \ar[-]{d} \\
A 
\end{tikzcd}
{\scriptsize\caption{\scriptsize The constituent structure of the set $\{\varnothing, \{\varnothing\}\}$.}}
\end{wrapfigure}
The purpose of the diagram is to show constituency relations rather than membership. The final step of the procedure above removes membership information which isn't implied by the constituency relations.

$A$ is a constituent of $B$ if $B$ can be reached from $A$ in the diagram by following a path of zero or more edges upward.

Every set contains the empty set as a constituent, so every set's symbol in the graph will be reachable from the empty set at the bottom via an upward path of edges.

Similarly, every set is a constituent of the full set, $S$, and is connected to $S$ by an upward path of edges.

As a simple example, we can consider the set $S=\{A, B\}$, with $B=\{A\}$ and $A=\varnothing=\{\}$. This set contains two elements, $A$ and $B$, and three constituents, $A$, $B$, and $S$ itself. $S$ is at the top of the diagram and $A$, the empty set, is at the bottom, with $B$ in the middle.

Despite the fact that $A$ is inside $S$ twice, once as an element of $S$ and once as an element of $B$, it appears in the diagram only once. Also, $A$ is an element of $S$, but there is no edge connecting $A$ directly to $S$. There is no way to tell, from the diagram, whether $A$ is an element of $S$ or not. The diagram only shows that $A$ is a constituent of $S$.

The constituent structure diagram doesn't show the full structure of $S$, including which constituents multiple levels below a set are also elements of it. That information has been removed because it is not part of the constituent structure{\footnote{The information which is missing from the constituent structure shown here is to be found in the category of sets, which knows the cardinality of each set shown, but nothing apart from that. The full information about the set, expressible using $\in$, has been separated into cardinality information, $\#$, and constituent information, $\lhd$.}}.

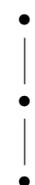
\begin{wrapfigure}{R}{3cm}
\centering
\begin{tikzcd}
\bullet \ar[-]{d}\\
\bullet \ar[-]{d} \\
\bullet 
\end{tikzcd}
\caption{\scriptsize The membership graph of $\{\{\varnothing\}\}$ is the same as its 
constituent structure graph.}
\end{wrapfigure}

Since we will often be interested in the constituent structures themselves rather than any specific sets with those structures, we may simply place a dot, $\bullet$, instead of a letter, at each position in the diagram.

In cases where we have a constituent structure, but no corresponding set, we can construct a set with that structure by starting with the empty set, identified with the dot at the bottom of the diagram, and assigning to each dot at any horizontal level the set containing as elements the sets denoted by the dots below it to which it is connected by a single edge. If two or more dots in the diagram are connected to the same collection of dots below, additional constituents can be selected from the dots farther down in the diagram, and then added as elements to the corresponding sets, to ensure that sets denoted by distinct dots in the diagram are distinct.

When it exists, the unique set whose membership graph is the same as a given constituent structure graph is the simplest set with that constituent structure. In later sections, we will develop algebraic procedures for constructing sets with specific constituent structures.

\section{Natural Numbers}

\subsection{Von Neumann's Construction}

Different ways of constructing the natural numbers from pure sets have been proposed. 
Von Neumann\cite{vonneumann} introduced the following system:
\begin{align}
0 \equiv& \{\}  \nonumber\\
1 \equiv& \{0\} \nonumber\\
2 \equiv& \{0, 1\} \\
3 \equiv& \{0, 1, 2\} \nonumber\\
4 \equiv& \{0, 1, 2, 3\} \nonumber\\
{} & \vdots \nonumber
\end{align}

\begin{wrapfigure}{R}{3cm}
\centering
\scriptsize
\begin{tikzcd}
3 \ar[-]{d} \\
2 \ar[-]{d} \\
1 \ar[-]{d} \\
0 
\end{tikzcd}
\caption{\scriptsize The constituent structure of von Neumann's construction of the natural number 3.}
\end{wrapfigure}

When we draw the constituent structure graph for the set identified with the number 3, the fact that it has more than one element is not visible in the graph. The sets corresponding to 2, 1 and 0 are all constituents of the set for 3, and, because of this, the fact that they are also elements of that set provides no constituent information.

The information about the set identified with 3 that can be expressed in terms of the constituency relation is exhaustively specified by $0\lhd 1 \lhd 2 \lhd 3$. It reveals that the set for 3 has at least one element, namely the set for 2, and at most three elements, since it has only three constituents apart from itself which could possibly be among its elements.

Somebody who knew the cardinality of each of these sets but didn't know their constituent structure would associate each of them with the number they are identified with; 3 has 3 elements and so on. From that point of view, this system of constructing the natural numbers successfully encodes the information in both the constituent structures of the sets and in their cardinalities.

\subsection{Zermelo's Construction}

Zermelo\cite{zermelo} used a different scheme for constructing the natural numbers:

\begin{align}
0 \equiv& \{\}  \nonumber\\
1 \equiv& \{0\} \nonumber\\
2 \equiv& \{1\} \\
3 \equiv& \{2\} \nonumber\\
4 \equiv& \{3\} \nonumber\\
{}  \vdots &\nonumber
\end{align}

\begin{wrapfigure}{R}{3cm}
\centering
\scriptsize
\begin{tikzcd}
3 \ar[-]{d} \\
2 \ar[-]{d} \\
1 \ar[-]{d} \\
0 
\end{tikzcd}
\caption{\scriptsize The constituent structure of Zermelo's construction of the natural number 3.}
\end{wrapfigure}
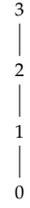

Apart from 0, each of these has only a single element, so somebody who only has access to cardinality information about a set would be unable to associate any of these sets with their corresponding number, apart from the empty set.

Considered in terms of constituent structure, however, each set is the simplest set with that structure. The edges in the constituent structure graph completely specify the membership relations between the sets. Zermelo's construction of the natural numbers is, in a sense, minimal, while von Neumann's is maximal, with each set containing as much as it possibly can.

From these observations, it is clear that, for each natural number, von Neumann's construction has the same constituent structure as Zermelo's construction of that number, despite the fact that the sets have different numbers of elements. The graph shared by both sets, for a natural number, $N$, always consists of a simple chain of $N+1$ constituents.

\subsection{Constituent Structure Isomorphisms}

When sets have the same number of elements, there is an invertible function which maps the elements of one set to the elements of the other set, establishing a one-to-one mapping between the sets. Those sets are ``the same'' in terms of cardinality.

In this case, Zermelo's sets and von Neumann's sets are ``the same'', but in a completely different way; they have the same constituent structures, but different cardinalities. There should therefore be an invertible map from one set to the other which preserves this structure in the same way that invertible functions between sets preserve the number of elements. It would not be a function between sets which puts their elements into one-to-one correspondence; instead it would map one constituent structure, depicted by the vertices and edges in the corresponding diagram, onto the other, sending vertices to vertices and edges to edges, in an invertible way.

\begin{wrapfigure}{R}{5cm}
\centering
\begin{tikzcd}
5_{{}_{Z}} \ar[-, ""{name=L1, below}]{d} \arrow[r] & 5_{{}_{V}}\ar[-, ""{name=R1, below}]{d}\\
4_{{}_{Z}} \ar[-, ""{name=L2, below}]{d} \arrow[r] & 4_{{}_{V}}\ar[-, ""{name=R2, below}]{d}\\
3_{{}_{Z}} \ar[-, ""{name=L3, below}]{d} \arrow[r] & 3_{{}_{V}}\ar[-, ""{name=R3, below}]{d}\\
2_{{}_{Z}} \ar[-, ""{name=L4, below}]{d} \arrow[r] & 2_{{}_{V}}\ar[-, ""{name=R4, below}]{d}\\
1_{{}_{Z}} \ar[-, ""{name=L5, below}]{d} \arrow[r] & 1_{{}_{V}}\ar[-, ""{name=R5, below}]{d}\\
0_{{}_{Z}}  \arrow[r]& {} 0_{{}_{V}}
\ar[shift left=2pt, from=L1, to=R1, ""]\ar[shift left=2pt, from=L2, to=R2, ""]\ar[shift left=2pt, from=L3, to=R3, ""]\ar[shift left=2pt, from=L4, to=R4, ""]\ar[shift left=2pt, from=L5, to=R5, ""]
\end{tikzcd}
\caption{\label{isomporphism}{\scriptsize{An isomorphic mapping, $M$, from the constituent structure of Zermelo's construction of the number 5 to that of von Neumann.}}}
\end{wrapfigure}
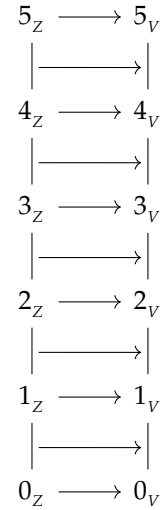

This mapping, considered to send one graph to the other, is a directed graph isomorphism. When the constituent structure is considered as a category, it's an invertible functor between the categories for the two sets. In the category of constituent structures, it's an isomorphism between objects.

Like a function, such a mapping sends elements of sets to elements of sets, but in this case, they are not all elements of the same set. It's the constituents of the two sets, rather than their elements, that are put into one-to-one correspondence by the isomorphism.

Since every element of a set is a constituent, and the isomorphism specifies an invertible mapping between the constituents of the source and the destination set, every element in the destination set must have something mapped to it from the source set, and every element of the source set must be mapped to something in the destination set.

Figure \ref{isomporphism} shows how the constituent structure of Zermelo's number 5 is mapped by an isomorphism, $M$, to the constituent structure of von Neumann's number 5, using subscripts of $Z$ and $V$ to denote Zermelo's construction and von Neumann's construction of each natural number.

The mapping is surjective, injective, and invertible, and also satisfies the condition:
\begin{equation}
a\lhd b \iff M(a)\lhd M(b).
\end{equation}

Because $5_{{}_{V}}$ has a larger cardinality than $5_{{}_{Z}}$, it has more {\em instances} of its constituents within it, including more elements. Multiple instances of a set are, however, literally the same thing. If $M(2_{{}_{Z}})=2_{{}_{V}}$, then every instance of $2_{{}_{V}}$ within $5_{{}_{V}}$ is the image of $2_{{}_{Z}}$, despite the fact that $5_{{}_{Z}}$ contains only one instance of $2_{{}_{Z}}$.

This does not break the rule that says that a mapping cannot send one thing to more than one thing. That rule exists to ensure that nothing is ever mapped to two different things. The multiple instances of $2_{{}_{V}}$ are not different things.

It would be a violation of that rule if two instances of the same set were to be mapped to two different destination sets. For example, in the inverse mapping, $M^{-1}{:}\ 5_{{}_{V}} \rightarrow 5_{{}_{Z}}$, a set such as $1_{{}_{V}}$ must be mapped to a single set, $M^{-1}(1_{{}_{V}})$. If two instances of $1_{{}_{V}}$, such as those which are elements of $2_{{}_{V}}$ and of $3_{{}_{V}}$ respectively, were to be mapped to different constituents of $5_{{}_{Z}}$, then it would truly be a case of one thing being mapped to two different things.

The mapping $M^{-1}{:}\ 5_{{}_{V}} \rightarrow 5_{{}_{Z}}$ is therefore an injective mapping from a five-element set to a one-element set, just as $M{:}\ 5_{{}_{Z}} \rightarrow 5_{{}_{V}}$ is a surjective mapping from a one-element set to a five-element set, which maps all of the constituents of $5_{{}_{Z}}$ to all of the constituents of $5_{{}_{V}}$ in a one-to-one way, despite there being 32 instances of the latter and only 6 instances of the former.

\subsection{A Constituent Structure Which is Not a Number}

The set produced by von Neumann's construction of a natural number such as 5 is very different from the corresponding set in Zermelo's scheme, but they have the same constituent structure.

It appears that the information that is lost when the constituent structure is extracted from each of those sets is the specification of which scheme was used for the construction, while the information kept by the constituent structure is the specification of the object constructed.

It is reasonable to enquire whether the different constructions of a single number will appear to be the same thing or different things within the constituent structure of a set that contains both of them.

\begin{wrapfigure}{R}{5cm}
\centering
\scriptsize
\begin{tikzcd}
{} &  \{2_{{}_{Z}}, 2_{{}_{V}}\} \ar[-]{ld} \ar[-]{rd} & {}\\
2_{{}_{Z}} \ar[-]{rd} & {} & 2_{{}_{V}} \ar[-]{ld} {} \\
{} & 1 \ar[-]{d} & {}\\
{} & 0 & {} 
\end{tikzcd}
\caption{\label{diamonddiagram}{\scriptsize{A set containing conflicting representations of the same natural number has a constituent structure which is not that of any natural number.}}}
\end{wrapfigure}

The simplest set containing two distinct representations of the same number is $\{2_{{}_{Z}}, 2_{{}_{V}}\}=\{\{1\}, \{0, 1\}\}$.
The constituent structure of this set is shown in figure \ref{diamonddiagram}. Neither of the sets, $2_{{}_{Z}}$ and $2_{{}_{V}}$, are constituents of each other, but $1=\{0\}$ and $0=\{\}$ are constituents of both. 

The fact that the set contains the two incompatible representations of the natural number 2 is clearly visible in the diagram of its constituent structure. The resulting graph is not isomorphic to the graph of any natural number.

This shows that the procedure of extracting the constituent structure from a set discards the information about which scheme was used to construct the natural numbers, as long as all occurrences of a natural number are encoded using the same set. Mixing one representation of a natural number with a distinct representation of the same number within a single set proves that something other than a consistent representation of the natural numbers was involved in the generation of that set.

The set $\{\{1\}, \{0, 1\}\}$ is the smallest, simplest set with a structure distinct from any natural number.
It can therefore be expected to play a role in the construction of other mathematical objects that are themselves distinct from natural numbers.

We will use the symbol $\Diamond$ to refer to this set, and call it the diamond set, in reference to its structure which distinguishes it from natural numbers, although in contexts where it plays an important specific role, we may refer to the same set in a different way to clearly state the role that it plays.

\section{Operations on Sets}

When we consider a set's constituents, the obvious ways in which two or more sets can be used in combination to specify another set are not the familiar operations of union and intersection, which are defined in terms of elements.

\subsection{Constituent Replacement}

The primary operation in which constituents play the central role is constituent replacement. We use the notation:
$$
x(y\rightarrow z)
$$
to denote the set that results from replacing every occurrence of $y$ in $x$ with $z$.

This operation can be performed for any three pure finite sets, $x$, $y$ and $z$. Any pure finite set can be expressed uniquely as a sequence of curly brackets, or braces, and commas\footnote{Equivalent representations of this type can be sorted and the first among them chosen as the unique one.}, and this operation is a simple substring substitution in that representation.

It has properties that are expressible using the $\lhd$ relation:
\begin{align}
y\lhd x &\implies z\lhd x(y\rightarrow z) \\
\lnot (y\lhd x) &\implies x(y\rightarrow z) = x \\
\lnot (y\lhd z) &\implies \lnot \left(y\lhd x(y\rightarrow z)\right)
\end{align}
where $\lnot$ is logical negation, indicating that the statement following it is false.

Every finite pure set contains the empty set, $\varnothing=\{\}$, as a constituent. In the special case in which the set being replaced is the empty set, we use the notation:
\begin{equation}
x(y) \  \equiv \ x(\{\} \rightarrow y).
\end{equation}

It has the properties:
\begin{align}
x(\{\}) = x \\
\{\}(x) = x \\
y\lhd x(y) \\
x(y)(y\rightarrow \{\}) = x \\
x(y)(z) = x(y(z)).
\end{align}

The properties above show that this binary operation is associative, invertible, has an identity element, $\{\}$, and that it constructs sets related by $\lhd$. $x(y)$ can be thought of as ``$x$ on top of $y$'', since its constituent structure diagram is $x$'s diagram on top of $y$'s diagram. 

It's also appropriate and helpful to think of $x(y)$ as ``$x$ after $y$'', since the procedure for constructing $x(y)$ starting from the empty set and successively enclosing sets within sets necessarily involves first constructing $y$, and then repeating the procedure for constructing $x$, but using $y$ as the starting point instead of the empty set.

\subsection{The Constituent Algebra}

The set $x(y)$ actually satisfies a stronger condition than $y\lhd x(y)$. There are no occurrences of the empty set within $x(y)$ other than those inside an occurrence of $y$. We can denote this by $x(y) \vdash y$. This notation can be thought of as graphically depicting $y$ as the horizontal line segment, $x$ as the vertical segment, and $x(y)$ as the entire $\vdash$ symbol, which is $x$ on top of $y$, displayed in the customary horizontal arrangement of symbols instead of vertically.

$\vdash$ can be defined as:
\begin{equation}
b\vdash a\  \equiv\  b(a \rightarrow \{\})(a)=b
\end{equation}
which means that every occurrence of the empty set within $b(a \rightarrow \{\})$ corresponds to an occurrence of $a$ within $b$.

When this condition, $b\vdash a$, is satisfied, there is some set, $c$, such that $b=c(a)$, which can be extracted from $b$ with the operation $b(a\rightarrow \{\})=c$, which ``removes $a$ from the bottom of $b$''. 

Similarly, when $b=c(a)$, the set $a$ can be obtained from $b$ and $c$ by ``removing $c$ from the top of $b$''. Consistency of notation suggests that we denote this as $(\{\}\leftarrow c)b=a$, and define it as:
\begin{equation}
(\{\}\leftarrow c)b\  \equiv
\begin{cases}
a & \exists a: b=c(a)\\
b & \rm{otherwise.} \\
\end{cases}
\end{equation}

This symmetry within the notation further suggests that we describe the relation between $x$ and $x(y)$ using the $\dashv$ symbol: $x\dashv x(y)$, which graphically depicts $x$ as the horizontal line segment, $y$ as the vertical one, and $x(y)$ as the entire $\dashv$ symbol, which consists of $x$ on top of $y$ when it is appropriately rotated to convert our left-to-right order of symbols into a top-to-bottom arrangement of constituent structures in a diagram, so:
\begin{equation}
c\dashv b\  \equiv\  \exists a: b=c(a)
\end{equation}
provides the definition of $\dashv$.

Note that the reversed symbol, $\dashv$, does not denote the same relation as $\vdash$ with the symbols in reverse order: $b\vdash a$ does not mean $a\dashv b$. The horizontal line segment in both cases points at a symbol whose constituent structure diagram is at the bottom, $\vdash$, or the top, $\dashv$, of that of the other symbol. So $b\vdash a$ can be expressed in words as ``$a$ is the bottom of $b$'', and $a\dashv b$ can be expressed as ``$a$ is the top of $b$''. 

The relations $\vdash$ and $\dashv$, like constituency, are partial orders on sets. They satisfy $\{\}\dashv x \vdash \{\}$ for all sets, $x$. The empty set is at the bottom of every set, $x$, because $x(\{\})=x$, so $x\vdash \{\}$, and at the same time, it's at the top of $x$ because adding it to the top of $x$ leaves $x$ unchanged: $\{\}(x)=x$, so $\{\}\dashv x$.

The notations for asserting that one set is at the top or bottom of another, and for removing one set from the top or bottom of another, lead to intelligible results, especially when we observe that associativity, $x(y(z))=x(y)(z)$, allows us to unambiguously use the expression $xyz$:
\begin{align}
xyz & \vdash z \\
x & \dashv xyz \\
xyz(z \rightarrow \{\}) & = xy \\
(\{\} \leftarrow x)xyz & = yz. 
\end{align}

The construction of a set from many others usually involves the addition of new sets to an expression such as $xyz$ on the left, to get a new set, $wxyz$, in which $w$ is added to the expression after $x$ is added. $wxyz$ can be read or thought of as $w$ after $x$ after $y$ after $z$, and the sequence of steps involved in the set's construction is encoded in the expression from right to left. 

Reading the expression $wxyz$ from right to left shows the order of its construction but conflicts with our conventions for reading and writing arithmetic expressions: $1+2+3$ is thought of as starting with 1 and then adding 2 and then 3. It is consistent with existing conventions for composition of functions, though: $f(g(h(x)))$ indicates that the functions should be evaluated starting with $h$ followed by $g$ and then $f$.

For finite pure sets in general or any other specific class of sets, we will refer to the expressions, operations, relations, and other sets that can be referenced, using the notation introduced here, as the constituent algebra of those sets.

\subsection{Extracting Constituents}

We can use the fact that $\lhd$ is a partial order to define an operation which selects the biggest constituents of a set.

For a set, $S$, we can say that a set, $m$, is maximal in $S$ if $S\neq m \lhd S$ and $m\lhd y \lhd S \Rightarrow y=m\ \rm{\sc or }\ y=S$. 

That is, a set is maximal in $S$ if there are exactly two distinct constituents of $S$ which contain it as a constituent, namely $S$ and itself.

We can write $\max{}^{\{\}}_\lhd(S)$ to denote the set which contains as elements all the maximal constituents of $S$.

In cases where there is only one maximal constituent we can call it the unique maximum and denote it by $\max{}_{\lhd}$. For example:
\begin{equation}
\max{}_{\lhd}(\{S\})=S.
\end{equation}

We don't need to resort to cardinality information and $\in$ to detect and extract the single maximal constituent from the set containing it. The set containing $y$ is $\{\varnothing\}(y) = \{\varnothing\}(\varnothing\rightarrow y)=\{y\}$, so we can extract the element $y$ from the set $\{y\}$ using $(\{\}\leftarrow \{\varnothing\})\{y\}=y$. We can detect whether there is only one element in a set, $x$, by seeing if $\{\varnothing\}\dashv x$.

In place of the intersection of two sets, the natural operation in this case is the selection of the Largest Common Constituents of two sets:
\begin{equation}
\rm{LCC}^{\{\}}(a, b) \equiv \max{}^{\{\}}_\lhd\{x\lhd a : x\lhd b\}.
\end{equation}

When there is only one of these, we can call it the Largest Common Constituent:
\begin{equation}
\rm{LCC}(a, b) \equiv \max{}_\lhd\{x\lhd a : x\lhd b\}.
\end{equation}

Finally, one can extract from a given set, $a$, the maximal constituents with another set, $b$, at the bottom:
\begin{equation}
a^{{}_{\{\}}}_{\vdash b} \equiv \max{}^{\{\}}_\lhd \{c \lhd a : c\vdash b\}
\end{equation}
and when there is only one, we can denote it by:
\begin{equation}
a_{\vdash b} \equiv \max{}_\lhd \{c \lhd a : c\vdash b\}
\end{equation}
and the constituents with $b$ at the top can be extracted in a similar way, but without using $\max{}^{\{\}}_\lhd$:
\begin{equation}
a^{{}_{\{\}}}_{b\dashv} \equiv \{c \lhd a : b\dashv c\}
\end{equation}
with the corresponding notation introduced to specify the unique constituent with $b$ at the top when there is only one: $a_{b\dashv} \equiv (\{\}\leftarrow \{\varnothing\})a^{{}_{\{\}}}_{b\dashv}$.

\section{Natural Number Arithmetic and the Guiding Principle}

The sets, $n_{{}_{Z}}$, that Zermelo identifies with natural numbers satisfy the following relation:
\begin{equation}
n_{{}_{Z}}(m_{{}_{Z}}) = n_{{}_{Z}} + m_{{}_{Z}}
\end{equation}
meaning that the set representing the sum of two natural numbers can be obtained by inserting the set representing one number into the set representing the other, replacing the empty set. 

The simplest non-trivial example is given by $1_{{}_{Z}}(1_{{}_{Z}})$. The set representing $1$ is $1_{{}_{Z}} = \{\varnothing\}$, which contains just the empty set, $\varnothing$, which represents zero. Replacing $\varnothing$ with $1_{{}_{Z}}$ yields $1_{{}_{Z}}(1_{{}_{Z}})=\{1_{{}_{Z}}\}=2_{{}_{Z}}$. The constituent replacement operation implements addition of natural numbers.

This establishes a simple correspondence between an elementary operation on sets and an elementary arithmetic operation on natural numbers.

The same relation is not true for von Neumann's construction. In that construction, each number is represented by a set with that number of elements. $2_{{}_{V}}$ has two elements, so regardless of which set replaces the single empty set inside the representation $1_{{}_{V}}=\{\varnothing\}$, the result of that replacement will always be a set with a single element, which can never be $2_{{}_{V}}$.

There is a different operation on set constituents which results in addition of natural numbers in von Neumann's representation. Replacing every constituent, $a$, of $n_{{}_{V}}$ with $a\cup m_{{}_{V}}$ produces the set $(n+m)_{{}_{V}}$. This can be expressed using the notation $n_{{}_{V}}(*\rightarrow *\cup m_{{}_{V}})=(n+m)_{{}_{V}}$.

Although the set operations which implement addition are different in the two cases, they produce the same effect on the constituent structures. For any two sets, $a$ and $b$, the diagram showing the constituent structure of $a(b)$ consists of the diagram for $a$ on top of the diagram for $b$, with the dot or symbol for the empty set at the bottom of $a$'s diagram identified with the dot or symbol for $b$ at the top of $b$'s diagram. The diagram for $a(*\rightarrow *\cup b)$ is exactly the same.

This is shown in table \ref{jointable} for two hypothetical sets, $a$ and $b$, with distinct constituent structures. 

\begin{table}
\caption{Set Operations Used to Join Constituent Structures Together}
\label{jointable}
\setlength{\tabcolsep}{5mm} 
\def\arraystretch{1.25} 
\centering
\scriptsize
\begin{tabular}{ccc}
\begin{tikzcd}
{} &  \bullet \ar[-]{ld} \ar[-]{rd} & {}\\
\bullet \ar[-]{rd} & {} & \bullet \ar[-]{ld} {} \\
{} & \bullet \ar[-]{d} & {}\\
{} & \bullet & {} 
\end{tikzcd} & \begin{tikzcd}
{} &  \bullet \ar[-]{ld} \ar[-]{rd} \ar[-]{d} & {}\\
\bullet \ar[-]{rd} & \bullet \ar[-]{d} & \bullet \ar[-]{ld} {} \\
{} & \bullet \ar[-]{d} & {}\\
{} & \bullet \ar[-]{d} & {}\\
{} & \bullet & {} 
\end{tikzcd} & \begin{tikzcd}
{} &  \bullet \ar[-]{ld} \ar[-]{rd} & {}\\
\bullet \ar[-]{rd} & {} & \bullet \ar[-]{ld} {} \\
{} & \bullet \ar[-]{d} & {}\\
{} &  \bullet \ar[-]{ld} \ar[-]{rd} \ar[-]{d} & {}\\
\bullet \ar[-]{rd} & \bullet \ar[-]{d} & \bullet \ar[-]{ld} {} \\
{} & \bullet \ar[-]{d} & {}\\
{} & \bullet \ar[-]{d} & {}\\
{} & \bullet & {} 
\end{tikzcd} \\
\scriptsize Constituent structure of $a$ & \scriptsize Constituent structure of $b$ & 
\scriptsize Constituent structure of $a(b)$ as well as $a(*\rightarrow *\cup b)$\\
\end{tabular}
\end{table}

It is worth observing that both the sets and the set operations are considerably more complicated in the case of von Neumann's construction of the natural numbers, due to the double burden of making the number visible in both the set's cardinality and its constituent structure.

The sets, $n_{{}_{Z}}$, used by Zermelo, on the other hand, are the simplest sets with the structure needed to represent the natural numbers. Their membership graphs, which specify everything about them, are identical with their constituent structure graphs, so they contain no information other than the structure of the natural number they represent.

The fact that addition of natural numbers is so simply implemented by a basic set operation, 
$n_{{}_{Z}}(m_{{}_{Z}}) = n_{{}_{Z}} + m_{{}_{Z}}$, when constituent structure alone is used to replicate the structure
of natural numbers, is a hint that this is the right way in general to construct and represent 
mathematical objects using sets\footnote{In the category whose objects are constituent structures and whose morphisms are constituent structure homomorphisms, this addition procedure corresponds to the sum of the two objects defined in terms of universal properties, while the product $x\star y$ is the structure obtained by replacing every edge in $x$ with a copy of $y$, identifying the top and bottom vertices in $y$'s diagram with the vertices at the top and the bottom of the edge in $x$. This operation implements multiplication of natural numbers, $n\star m = n\times m$, and results in a set only when a procedure for constructing a set from the resulting constituent structure is specified. The simplest set with that structure is the natural choice, which makes sums and products of Zermelo's natural numbers coincide with sums and products of constituent structures of sets in general.}.

This hint prompts us to observe that natural numbers are the only structures which it is possible to encode in the cardinalities of finite sets. Anything more complicated, such as a negative integer or an ordered pair, will need to be encoded in a set's constituent structure, since the cardinality is always a natural number.

In fact, mathematical objects in general have structures, while sets have constituent structures, elements and cardinality. When we construct an object from sets and subsequently abstract from the underlying set to get the object as an object rather than as a set, the resulting object has structure but no elements or cardinality. The constituent structure of a set is the structure which survives within the constructed object when that abstraction occurs.

This gives us a guiding principle when constructing mathematical objects in general from sets:

\begin{itemize}
\item The structure of the object must be encoded solely in the set's constituent structure. 
\end{itemize}

Everything else will be lost when we forget that the underlying object is a set.

\section{Ordered Pairs and Tuples}
\label{orderedpairs}

\subsection{The Standard Construction of an Ordered Pair}

First introduced by Kuratowski\cite{kuratowski}, the definition of an ordered pair of sets as:
\begin{equation}
(a, b) = \{\{a\}, \{a, b\}\}
\end{equation}
is universally accepted today for good reason.

It is extremely simple, intuitive, and it specifies the order of $a$ and $b$ successfully for all sets. Its structure clearly encodes the appropriate concept: First $a$, then $b$.

\begin{wrapfigure}{R}{5cm}
\centering
\scriptsize
\begin{tikzcd}
{} &  \{\{a\}, \{a, b\}\} \ar[-]{ld} \ar[-]{rd} & {}\\
\{a\} \ar[-]{d} & {} & \{a, b\} \ar[-]{lld} \ar[-]{d} \\
a  \ar[-]{dr} & {} & b \ar[-]{dl} \\
{} & \rm{LCC}(a, b) & {} 
\end{tikzcd}
\caption{\label{kuratowskisucceeds}{\scriptsize{When neither $a$ nor $b$ are constituents of each other, their order and structures can be retrieved from the Kuratowski pair's structure.}}}
\end{wrapfigure}

We bring new requirements, though:

\begin{itemize}
    \item The set $(a, b)$ should have a constituent structure which contains within it the constituent structures of $a$ and $b$, in the correct order and separately retrievable.
    \item Given the set $(a, b)$, it should be clear from its constituent structure diagram that its structure is that of an ordered pair of two sets.
\end{itemize}

So we are considering the case when we can't see the elements of $(a, b)$, only the constituents. We have access to the information regarding $(a, b)$ expressible in terms of $\lhd$, which we can display in a diagram, and from that information, we need to be able to reconstruct the diagrams for $a$ and $b$ and in addition determine which of the two is first and which is second in the pair.

\begin{wrapfigure}{R}{5cm}
\centering
\scriptsize
\begin{tikzcd}
{} &  \{\{a\}, \{a, b\}\} \ar[-]{ld} \ar[-]{rd} & {}\\
\{a\} \ar[-]{rd} & {} & \{a, b\} \ar[-]{ld} {} \\
{} & a \ar[-]{d} & {}\\
{} & \bullet \ar[-]{d}  & {} \\
{} & \vdots & {}
\end{tikzcd}
\caption{{\scriptsize{When $b\lhd a$, $a$'s structure can be recovered, but $b$ could be any constituent of $a$ other than $a$.}}}
\end{wrapfigure}

One indication that Kuratowski's definition might not be sufficient for us can be seen from the fact that, in order to tell which of the two sets in $\{\{a\}, \{a, b\}\}$ contains the first set in the pair, we use the cardinalities of the sets, choosing the set with one element. Cardinality information is exactly what we need to avoid using.

The shape of the constituent structure diagram for $\{\{a\}, \{a, b\}\}$ will depend on the details of $a$ and $b$. In the best case, neither entry in the pair will be a constituent of the other, and the resulting constituent structure is shown in figure \ref{kuratowskisucceeds}. The unknown internal structures of $a$ and $b$ are omitted from the diagram, and it has been assumed for simplicity that they have a unique largest common constituent, whose structure is also omitted.

In this case, it's clear from the fact that $a$ is within both elements of the pair's set, while $b$ is in just one, that $a$ is the first entry and $b$ is the second.

It's also possible to recover the constituent structure of $a$ from the graph; it includes $a$ and everything which can be reached from $a$ along a downward path of edges. $b$'s structure can be recovered in the same way.

Next, consider the case when $b\lhd a$. $a$ can be identified as the set contained in both elements of the set representing the ordered pair, but $b$ is among the constituents of $a$, with no way to tell from the diagram which of those constituents of $a$ corresponds to $b$.

\begin{wrapfigure}{R}{5cm}
\centering
\scriptsize
\begin{tikzcd}
\{\{a\}, \{a, b\}\} \ar[-]{d} \\
\{a, b\} \ar[-]{d} \\
b \ar[-]{d} \\
\bullet \ar[-]{d}\\
\vdots
\end{tikzcd}
\caption{{\scriptsize{When $\{a\}\lhd b$, only $b$'s structure can be recovered, and the diagram's structure does not indicate that it encodes an ordered pair.}}}
\end{wrapfigure}

There is only one case in which the constituent of $a$ which is equal to $b$ can be deduced, specifically if $a$ has only one constituent apart from itself, which would have to be the empty set. So the Kuratowski ordered pair has a unique pair of sets, namely $b=\varnothing$ and $a=\{\varnothing\}$, which it is able to represent as arranged in a specific order, $\left(\{\varnothing\}, \varnothing\right)$. 

This uniquely representable ordering of the two simplest possible sets, $0=\{\}$ and $1=\{0\}$, is encoded in the set $\left\{\{1\}, \{1, 0\}\right\}=(1, 0)$, which is the $\Diamond$ set encountered previously, the simplest set which does not have the structure of a natural number.

The Kuratowski pair also fails to encode both entries in the pair when $\{a\}\lhd b$. $b$ is then recoverable from the constituent structure since it is the vertex in the third horizontal layer, but $a$ is not specified by the diagram, and it is also not evident from the diagram that it encodes an ordered pair.

\subsection{Construction of an Ordered Pair Which Respects Constituent Structure}

From these examples it is clear that a different way of encoding ordered pairs in a set is necessary in order to encode arbitrary constituent structures in a specific order in another constituent structure in a way that allows them to be retrieved.

One difficult case which indicates the general form that the solution must have is the case when one of the sets is several layers deep and the other is the empty set. The empty set is always at the very bottom of the graph, with a single upward connection to $\{\varnothing\}$. 

The only way for the diagram itself to point at a specific vertex is by connecting an edge to it, but we need to be able to see that edge and distinguish it from the edges which form part of the structures of the entries in the pair.

There is no way for the shape of the diagram in higher layers to connect such a distinguished edge to the empty set, since all paths to the empty set go through $\{\varnothing\}$.

From this, we can conclude that a vertex higher up in the graph will need to act as a proxy for the empty set. If the first entry in the pair, $a$, is the empty set, $\varnothing$, then there would need to be a proxy vertex, $v_\varnothing$, which is not a constituent of the second entry, $b$, so that an edge from the highest levels of the graph can connect to $v_\varnothing$ without travelling through any part of $b$, to indicate that $v_\varnothing$ is the proxy for the first entry in the ordered pair. 

These considerations lead to an inevitable conclusion: $a$ must be placed on top of something, $v_a$, resulting in the set $a(v_a)$, and $b$ must be placed on top of something else, $v_b$, resulting in $b(v_b)$, in a way that makes it impossible for $a(v_a)$ to contain $v_b$ or for $b(v_b)$ to contain $v_a$.

This requires the use of two sets, $v_a$ and $v_b$, whose constituent structures are incompatible with the possibility of either one being a constituent of a set with the other at the bottom.

We have an example of a set, denoted by $\Diamond$, with a structure that can never be isomorphic to a natural number. So if $a$ is placed on top of that set, and $b$ is placed on top of a natural number, such as $3_{{}_{Z}}$, then neither of the resulting sets, $a(\Diamond)$ and $b(3_{{}_{Z}})$, can contain the other.

This scheme only works for ordered pairs: With an ordered triple, $(a, b, c)$, $c$ can't be placed on top of a natural number, since it could then contain or be contained in $b(3_{{}_{Z}})$, and it can't be placed on top of a diamond, since it could then contain or be contained in $a(\Diamond)$.

The general solution, which, in the case of an ordered tuple, $(x_0, x_1, \cdots, x_k)$, generates corresponding sets which are guaranteed not to contain each other, is to place each set on top of $\Diamond$ on top of a distinct natural number, generating the set $x_n(\Diamond(n))$ for the $n^{\rm{th}}$ entry.

In this context, $\Diamond$ takes on the role of a position indicator, and its role is more clearly comprehensible when it is denoted by the word Position.

With this notation, the ordered pair with $a$ in the first position and $b$ in the second position is successfully encoded in the constituent structure of the set:
\begin{equation}
(a, b) = \{a({\rm Position}(0)), b({\rm Position}(1))\}.
\end{equation}

The constituent structures of the two sets used here are shown in figure \ref{orderedpairwithpositions}.

{\begin{figure}
\centering
\scriptsize
\begin{tikzcd}
{} &  {}  & {}   & (a, b) \ar[-, gray, dashed]{lldd} \ar[-, gray, dashed]{rrd} & {} & {} & {}   \\
{} &  {}  & {}   & {}  & {} & b({\rm Position(1)}) \ar[-]{d} & {}   \\
{} &  a({\rm Position}(0)) \ar[-]{d}  & {}   & {}  & {} & {\rm Position(1)} \ar[-]{ld} \ar[-]{rd}& {}   \\
{} &  {\rm Position(0)} \ar[-]{ld} \ar[-]{rd} & {}   & {}  & \bullet \ar[-]{rd}& {} & \bullet  \ar[-]{ld} \\
\bullet \ar[-]{rd} & {} & \bullet \ar[-]{ld} {}   & {}  & {} & \bullet \ar[-]{d} & {}   \\
{} & \bullet \ar[-]{d} & {}                       & {}  & {} & 1 \ar[-]{d} & {}   \\
{} & 0 & {}                       & {}  & {} & \bullet  & {}   \\
\end{tikzcd}
\caption{The sets $a({\rm Position}(0))$ and $b({\rm Position}(1))$ encode the constituent structures of $a$ and $b$ in their specified positions with no possibility that one could be a constituent of the other. The constituent structure diagrams of the two sets are shown side-by-side for clarity, rather than as constituents of a single set containing them both. Grey dashed lines indicate edges that would connect the two sets to a set that did contain them as elements. The constituent structure diagram for such a set would fuse the bottom three vertices of the two diagrams shown here together.
\label{orderedpairwithpositions}}
\end{figure}
}

\subsection{Extracting and Addressing Sets Inside Tuples of Tuples}

A set, $T_0$, containing several distinct sets, $s_0, s_1, \cdots, s_k$, in their respective positions in that order, is:
\begin{equation}
T_0 = \{s_0({\rm Position}(0)), s_1({\rm Position}(1)), \cdots, s_k({\rm Position}(k))\}
\end{equation}
which successfully encodes the constituent structures of those sets within its own, along with their order,
which provides a construction of the (k+1)-tuple $(s_0, s_1, \cdots, s_k)$ as a set.

That tuple may then be included within another tuple, $T_1$. For example, it might be put at the second position of an ordered pair after a different set, $U$:
\begin{equation}
T_1 = \{U({\rm Position}(0)), T_0({\rm Position}(1))\}.
\end{equation}

The set $T_0({\rm Position}(1))$ is $T_0$ with all of the empty sets inside it replaced by ${\rm Position}(1)$.
Those empty sets inside $T_0$ are deep within its elements, such as $s_0({\rm Position}(0))$. Replacing the empty sets inside $s_0({\rm Position}(0))$ results in $s_0({\rm Position}(0){\rm Position}(1))$, which is therefore an element of $T_0({\rm Position}(1))$ as well as a constituent of $T_1$.

If the tuple, $T_1$, is subsequently positioned within another tuple, which in turn is positioned within another, and so on, the resulting $(n+1)$-dimensional tuple, $T$, will contain the constituent:
$$
s_0({\rm Position}(p_0){\rm Position}(p_1)\cdots{\rm Position}(p_n))
$$
where $p_j$ is the position in the $j^{\rm{th}}$ tuple of the previous tuple for $j>0$.

If we introduce the notation ${\rm Position}(p_0, p_1, \cdots, p_n)$ to refer to the multi-dimensional position set ${\rm Position}(p_0)\cdots{\rm Position}(p_n)$, then we can assert that $s_0$ is at that $(n+1)$-dimensional position within the $(n+1)$-dimensional tuple, $T$, with the expression:
\begin{equation}
\label{objectatposition}
s_0({\rm Position}(p_0, p_1, \cdots, p_n)) \lhd T.
\end{equation}

Every constituent of $s_0$ satisfies the same condition, so this actually asserts that $s_0$ is a constituent of the object at that position. This means that we can assert that $T$ contains {\em something} at that position with:
\begin{equation}
{\rm Position}(p_0, p_1, \cdots, p_n) \lhd T.
\end{equation}

To extract the full set at that position, we can take the maximal constituent of $T$ with that position at the bottom:
\begin{equation}
\label{setatmultiposition}
T_{{\rm\vdash  Position}(p_0, p_1, \cdots, p_n)} = s_0({\rm Position}(p_0, p_1, \cdots, p_n)).
\end{equation}

Equation \ref{setatmultiposition} also allows us to assert that $s_0$ is the full set at that position and not just a constituent of that set.

This gives us $s_0$ on top of its multi-position set. To retrieve $s_0$ in its original form, we can remove the position set from the bottom of the constituent of $T$:
\begin{equation}
s_0= T_{{\rm \vdash Position}(p_0, p_1, \cdots, p_n)}({\rm Position}(p_0, p_1, \cdots, p_n) \rightarrow \{\}).
\end{equation}

\subsection{Observations}

The construction of ordered pairs and tuples given here has a number of desirable properties, one of which is that it is essentially the only solution. Minor adjustments, such as starting with position 1 instead of 0, or using a more complicated not-a-number object in place of $\Diamond$, are possible, but any solution which works for all sets will involve placing the entries in the tuple on top of sets with non-isomorphic structures, of which the ones given here are the simplest. This is dictated by the requirement that the order and constituent structures of the entries should be retrievable from the constituent structure of the tuple.

The result of adherence to that requirement is a positional representation which closely matches our existing concept of an ordered pair: It's a structure containing two things, one of which is in the first position and the second of which is in the second position\footnote{One possible variation of this method is to use the entries themselves instead of numbers to specify the order. For example, when $a$ is not empty, $\{\{a\}, b\{\varnothing, a\}\}$ places $b$ after $a$ rather than in position 1, and $b$ can be retrieved by querying for the constituent with $\{\varnothing, a\}$ at the bottom, giving this the structure of a linked list rather than an ordered pair.}.

The other benefit is that the expressions needed to construct the ordered pairs and even multi-dimensional tuples, to address the positions within them and the objects at those positions, and to extract the objects from them, are extremely simple and also closely match our own conceptual understanding of what is involved in each case.

Like the case of addition for natural numbers, when the structure of each object is replicated solely by the constituent structure of the underlying set, the operations and relations natural to the objects represented coincide with simple relations and operations on sets.

It is worth bearing in mind that we were constrained to this form of solution because it was necessary to be able to encode and retrieve arbitrary sets. For sets of a specific type, such as those which encode natural numbers in Zermelo's encoding, easier methods of encoding ordered pairs of those sets may exist.

In fact, the construction of multi-dimensional tuples given here specifies a way to encode finite sequences of natural numbers. The multi-dimensional position set, ${\rm Position}(p_0, p_1, \cdots, p_n)$, is essentially a finite sequence of natural numbers with a non-number object, $\Diamond$, situated between each and the next: ${\rm Position}(p_0, p_1, \cdots, p_n)=\Diamond p_0\Diamond p_1\cdots\Diamond p_n$.

In this context, the $\Diamond$ set plays the role of a comma, delimiting the numbers which specify a sequence of coordinates.

That set itself is a unique representation of two specific objects in a specific order. It's the Kuratowski ordering of 1 followed by 0, and these are the only sets already ordered by constituency which can unambiguously be given a Kuratowski ordering within a constituent structure, and it can only be this order: 1, then 0.

This reverse ordering of the two simplest natural numbers, which appears in the simplest set whose structure isn't that of any natural number, plays a central role in Part Two of this series, where we consider numbers whose structure is more complicated than a simple count.

\section{Top and Bottom Structures}

\subsection{Encoding Part of a Constituent Structure}

The set operation $xy=x(\{\}\rightarrow y)$ allows us to combine two constituent structures to obtain another which consists of $x$'s structure on top of $y$'s, by fusing them together at a single vertex. Not all constituent structures can be generated by doing this, though, and that set operation can't construct all sets. 

Figure \ref{topandbottomstructures} shows an example of parts which might appear at the top and at the bottom of a constituent structure diagram. Neither of these are constituent structures of actual sets, but the vertices labelled $x$, $y$, and $z$ designate distinct sets with specific constituent structures, and the diagram obtained by identifying those with the vertices labelled $a$, $b$, and $c$ in the top diagram is a constituent structure of an actual set. Following the convention in graph theory, we can refer to these unconnected vertices as terminals.

{\begin{figure}
\centering
\scriptsize
\begin{tikzcd}
        &{}      &{} \bullet \ar[-]{dl} \ar[-]{dr} &{}      \\
        &{} \bullet \ar[-]{dl} \ar[-]{dr}    &{}     &{}  \bullet \ar[-]{d} \\
a &{}      &{} b  &c    \\
{}         &{}             &y      &z      \\
x   &{}  &\bullet \ar[-]{u}\ar[-]{ur}  &{}  \\
{}         &{}  \bullet \ar[-]{ul}\ar[-]{ur}    &{}&{}\\
{}         &{} \bullet \ar[-]{u} &{}         &{}  \\
\end{tikzcd}
\caption{Structures which may appear at the top or bottom of a constituent structure diagram. A single diagram can be constructed by identifying the vertices labelled $a$, $b$, and $c$ in the top diagram with those labelled $x$, $y$, and $z$ in the bottom diagram. 
\label{topandbottomstructures}}
\end{figure}
}

When $x$ and $y$ are sets which aren't constituents of each other, the Kuratowski ordered pair $(x, y)_K=\{\{x\}, \{x, y\}\}$ always has the same shape at the top of its constituent structure diagram, which can be seen in figure \ref{kuratowskisucceeds}, while the shape of the bottom depends on $x$ and $y$. The set $d=ab$ has $a$ at the top, which we can express as $a\dashv d$, and $(x, y)_K$ has something specific at the top as well, which the notation we have introduced so far isn't able to express.

We will need to introduce some new notation into the constituent algebra to allow it to describe and fuse together structures with more than one terminal appearing at the top or bottom of their constituent structure diagrams.

The techniques introduced in the previous section can be used again here to ensure that the operations we describe can be performed with arbitrary sets without losing any constituent information.

For any structure which could appear at the top of a constituent structure diagram, we can specify a set which has that structure at the top by attaching $\Diamond n=\rm{Position}(n)$ to the $n^{\rm th}$ terminal\footnote{Technically, a structure could have so many dots within it at the same horizontal level that, even after the sets $\Diamond n$ have been attached to its terminals, there are more dots in the diagram than distinct combinations of constituents available to serve as elements of the corresponding sets, each of which needs to have a unique combination of elements in order to be a distinct set represented by a distinct dot. In those cases, there is a natural number, $m$, such that the sets $\Diamond (n+m)$ can be attached to the terminals to achieve the desired result.}, as shown in figure \ref{topstructure}.

{\begin{figure}
\centering
\scriptsize
\begin{tikzcd}
        &{}      &{} T \ar[-]{dl} \ar[-]{dr} &{}      \\
        &{} \bullet \ar[-]{dl} \ar[-]{dr}    &{}     &{}  \bullet \ar[-]{d} \\
\Diamond 0 \ar[-]{dddrr} &{}      &{} \Diamond 1 \ar[-]{ddr} & \Diamond 2  \ar[-]{d}  \\
{} &{}      &{}  & \bullet \ar[-]{d}    \\
 &{}       &{}  & \bullet \ar[-]{dl}   \\
{} &{}      &{} \bullet \ar[-]{d} & {}    \\
{} &{}      &{}  \bullet    &{}    & {}   \\
\end{tikzcd}
\caption{The constituent structure of the simplest set, $T$, with the structure shown in figure \ref{topandbottomstructures} at its top. The sets which are fused to the terminals $a$, $b$, and $c$, are not constituents of one another or of any set $\rm{Position}(n)=\Diamond n$ with $n>2$. The structure of each $\Diamond$ set has been shrunk to a single vertex in each instance for clarity.
\label{topstructure}}
\end{figure}
}

The set, $T$, shown in the figure has distinct labelled and ordered vertices corresponding to the terminals $a$, $b$ and $c$ from figure \ref{topandbottomstructures}, consisting of the sets $\Diamond n=\rm{Position}(n)$, where $n$ is a natural number below 3. $T$ is the most natural encoding within a set of that partial constituent structure.

$T$ has the properties that it has $m>0$ constituents which are of the form $\Diamond n$, with $n<m$, and its diagram has no paths connecting $T$ to $0$ which don't go through any of those constituents: 
$$
x\lhd T \implies x \lhd \Diamond n \mbox{\sc\small\ \ or\ \ } \Diamond n \lhd x.
$$

We will call a set with these properties a top structure, indicating that it encodes the specification of the top part of a constituent structure diagram, and we will call the constituents $\rm{Position}(n)$ its terminals.

This construction allows us to identify the set which naturally encodes the top of the Kuratowski pair's diagram as $(\rm{Position}(0), \rm{Position}(1))_K$, which is conceptually comprehensible as a specification of a Kuratowski ordered pair whose first entry is whatever is at position zero and whose second entry is whatever is at position one.

An analogous procedure identifies a set, $B$, which naturally encodes any specified bottom part of a constituent structure. The set $\{x, y, z\}$ has the lower structure shown in figure \ref{topandbottomstructures} at its bottom, but $x$, $y$, and $z$ have no specific ordering or labelling within its constituent structure diagram which allows them to be matched with the corresponding terminal in $T$.

The set $\{x\ {\rm Position}(0), y\ {\rm Position}(1), z\ {\rm Position}(2)\}$ naturally encodes the ordered triple $(x, y, z)$ and provides all of the needed information, but it doesn't have the specified structure at the bottom of its constituent structure diagram, and if we are given that set, we have no way of knowing that it is intended to specify the bottom of a constituent structure rather than an ordered triple of sets.

We also cannot use the set $\{{\rm Position}(0) x, {\rm Position}(1) y, {\rm Position}(2) z\}=\{\Diamond x, \Diamond1 y, \Diamond 2z\}$, which has the right structure at the bottom, because the sets $x$, $y$, and $z$ could be of the form $x=2g,y=1h, z=k$, in which case the resulting set, $\{\Diamond2g, \Diamond2h, \Diamond2k\}$, wouldn't encode the order of the terminals.

{\begin{figure}
\centering
\scriptsize
\begin{tikzcd}
        &{}      &{} B \ar[-]{lldddd} \ar[-]{dd} \ar[-]{dr}&{}       \\
        &{}      &{} &2\Diamond\ z  \ar[-]{d}     \\
  {}  &{}    &{} 1\Diamond\ y    \ar[-]{d}     &{}  1\Diamond\ z \ar[-]{d}  \\
{} &{}      &{} \Diamond\ y \ar[-]{d}& \Diamond\ z  \ar[-]{d}  \\
\Diamond\ x    \ar[-]{d}     &{}             &y      &z      \\
x   &{}  &\bullet \ar[-]{u}\ar[-]{ur}  &{}  \\
{}         &{}  \bullet \ar[-]{ul}\ar[-]{ur}    &{}&{}\\
{}         &{} \bullet \ar[-]{u} &{}         &{}  \\
\end{tikzcd}
\caption{The set, $B$, which naturally encodes the bottom structure of figure \ref{topandbottomstructures}, using $n\Diamond s$ to specify that the set $s$ should be placed in the $n^{\rm th}$ position in a top structure. The internal structure of each $\Diamond$ set has been omitted for clarity.
\label{bottomstructure}}
\end{figure}
}

The set that works is $B=\{(0){\rm Position}\ x, (1){\rm Position}\ y, (2){\rm Position}\ z\}=\{\Diamond x, 1\Diamond y, 2\Diamond z\}$. The structure of this set is shown in figure \ref{bottomstructure}.

The terminals in $B$ can be identified and matched with the terminals in $T$. A constituent, $C_n$, of $B$  matches a terminal $T_n=\Diamond n$ of $T$ if $C_n$ is maximal in $B$ and $T_n\Diamond \dashv \Diamond C_n$.

The first condition ensures that the number at the top of $C_n$ is not just a part of a larger number. The second condition, $T_n\Diamond \dashv \Diamond C_n$, ensures that the number at the top of $C_n$ is neither larger nor smaller than the number at the bottom of $T_n$.

The corresponding terminal, $B_n$, can be extracted from $C_n$ and the matching terminal, $T_n$, using:
\begin{equation}
B_n = (\{\}\leftarrow T_n\Diamond)(\Diamond C_n)
\end{equation}
or more simply:
\begin{equation}
B_n = (\{\}\leftarrow n\Diamond) C_n.
\end{equation}

We will say that a set, $B$, is a bottom structure if, for some natural number $m>0$, it has $m$ maximal constituents, $C_n$, with distinct values of $n<m$, such that $C_n=n\Diamond x_n$ for some set $x_n$, and has no path in its constituent structure diagram that connects it to zero without going through one of these constituents:
$$
x\lhd B \implies x \lhd C_n \mbox{\sc\small\ \ or\ \ } C_n \lhd x.
$$

\subsection{Fusing Structures Together}
\label{fusing}

Given a top structure, $T$, and a bottom structure, $B$, with the same number of terminals, such as those shown in figures \ref{topstructure} and \ref{bottomstructure}, there is a specific structure which results from connecting the terminals of $T$ to the corresponding terminals of $B$.

We can introduce the operation which accomplishes this into the constituent algebra by specifying the set operations that need to be performed on $T$ and $B$ to produce the set with that constituent structure.

The obvious way to do this is simply to use constituent replacement to replace each of $T$'s terminals with the corresponding terminal in $B$:
$$
T(\Diamond \rightarrow B_0)(\Diamond1 \rightarrow B_1)\cdots (\Diamond m\rightarrow B_m)
$$
which we can write in a more compact way with the notation:
$$
T(\Diamond n \rightarrow B_n)^m_{n=0}
$$
but when replacing constituents successively for each terminal, there is a risk that the later replacement operations could change the sets which were previously inserted.

For example, $x$ might contain $\Diamond 1$ as a constituent, so the set:
$$
T(\Diamond \rightarrow x)(\Diamond 1 \rightarrow y)(\Diamond 2 \rightarrow z)
$$
contains the constituent $x(\Diamond 1 \rightarrow y)(\Diamond 2 \rightarrow z)$, which is not necessarily the same set as $x$.

This risk is avoided by inserting the terminals of $B$ on top of the terminals of $T$:
$$
T(\Diamond n \rightarrow B_n\Diamond n)^m_{n=0}
$$
after which the terminals of $T$ can be removed:
$$
T(\Diamond n \rightarrow B_n\Diamond n)^m_{n=0}(\Diamond n \rightarrow \Diamond (n-1))^1_{n=m}(\Diamond \rightarrow \{\})
$$
where the replacements $\Diamond n \rightarrow \Diamond (n-1)$ ensure that $B_n$ is never a constituent of a set which has a replacement operation performed on it.

As these replacements occur, the constituent $B_n\Diamond n$ changes to $B_n\Diamond (n-1)$, $B_n\Diamond (n-2)$ and so on, until the final replacement which changes $B_n\Diamond$ to $B_n$, after which every terminal of $B$ is a constituent of the resulting set.

\subsection{Notation}

We can denote the operation that fuses a top structure to a bottom structure with a double underline:
\begin{equation}
\label{fusenotation}
\underline{\underline{TB}} \ \equiv \ \ T(\Diamond n \rightarrow B_n\Diamond n)^m_{n=0}(\Diamond n \rightarrow \Diamond (n-1))^1_{n=m}(\Diamond \rightarrow \{\})
\end{equation}
where $m$ is the number of terminals in $B$ and also in $T$. This notation indicates that more than one connection is being made between $T$ and $B$, and it also groups the symbols together so that in an expression such as $w\underline{\underline{xy}}z$, the expression $\underline{\underline{xy}}$ denotes a set with a fully specified constituent structure that has no terminals, like $w$ and $z$.

The set that results from this operation has one of the structures at the top and the other at the bottom. We can denote this relation between a set and a top or bottom structure using the notation:
\begin{equation}
t\Dashv x \ \equiv\ \exists b : x= \underline{\underline{tb}}
\end{equation}
and:
\begin{equation}
x\vDash b \ \equiv\ \exists t : x= \underline{\underline{tb}}
\end{equation}
where the two lines in the double turnstile symbol indicate that there is more than one connection between the bottom and the top structure in each case.

So in the example of the Kuratowski ordered pair, we can use the symbol $K$ to denote the structure at the top of every such pair, $K=(\rm{Position}(0), \rm{Position}(1))_K$, and express the fact that a given pair has that structure at the top with:
\begin{equation}
K\Dashv (x, y)_K.
\end{equation}

Finally, we can call a set which is both a top structure and a bottom structure, and whose top structure terminals are distinct from its bottom structure terminals, a middle structure. The fusing operation defined above can combine a middle structure with a top structure to produce a different top structure, or with a bottom structure to produce a new bottom structure. This can be naturally included in an expression which constructs a set by using the notation $v\underline{\underline{wxy}}z$, which is unambiguous because the operation $\underline{\underline{xy}}$ is associative.

\subsection{Composing Middle Structures}

The ordered pairs and tuples constructed in section \ref{orderedpairs} satisfy the definition of a top structure: $(a, b, c)=\{a\Diamond, b\Diamond1, c\Diamond2\}$, while a set of the form $\{\Diamond x, 1\Diamond y, 2\Diamond z\}$ is a bottom structure.

The set: $$m_1=\{\Diamond x\Diamond,1\Diamond y\Diamond1,2\Diamond z\Diamond2\}$$ has both of these forms, and is a middle structure. We can fuse it to another middle structure, such as $m_2=\{\Diamond a\Diamond,1\Diamond b\Diamond1,2\Diamond c\Diamond2\}$, to get a third middle structure:
\begin{equation}
\underline{\underline{m_1m_2}} = \{\Diamond xa\Diamond,1\Diamond yb\Diamond1,2\Diamond zc\Diamond2\}
\end{equation}
or we could fuse them in the reverse order to get a different set:
\begin{equation}
\underline{\underline{m_2m_1}} = \{\Diamond ax\Diamond,1\Diamond by\Diamond1,2\Diamond cz\Diamond2\}.
\end{equation}

Middle structures with a given number of terminals form an algebra with an associative, non-commutative, binary operation, similar to the binary operation on sets which replaces the occurrences of the empty set in one set with the other set.

There is an identity element for this operation on middle structures: 
\begin{equation}
I_{\mbox{\tiny{$M$}}} =  \{\Diamond \Diamond,1\Diamond \Diamond1, \cdots, n\Diamond \Diamond n\}
\end{equation}
which satisfies:
\begin{equation}
\underline{\underline{I_{\mbox{\tiny{$M$}}}x}} = \underline{\underline{xI_{\mbox{\tiny{$M$}}}}} = x
\end{equation}
for any middle structure, $x$, with $n$ terminals.

A middle structure whose top terminals have natural numbers which differ from those of the bottom terminals to which they are attached permutes the terminals of other middle structures:
\begin{align}
p & = \{\Diamond\Diamond1,1\Diamond \Diamond2, 2\Diamond \Diamond \} \nonumber\\
\underline{\underline{pm_1}} &= \{\Diamond y\Diamond1,1\Diamond z\Diamond2,2\Diamond x\Diamond\} \nonumber\\
\underline{\underline{m_1p}} &= \{1\Diamond y\Diamond2,2\Diamond z\Diamond,\Diamond x\Diamond1\} \\
\underline{\underline{pm_1p}} &= \{\Diamond y\Diamond2, 1\Diamond z\Diamond,2\Diamond x\Diamond1\} \nonumber\\
\underline{\underline{pm_1pp}} &= \{\Diamond y\Diamond, 1\Diamond z\Diamond1,2\Diamond x\Diamond2\} \nonumber
\end{align}
where the fact that the permutation $p$ is a cycle of length 3 shows in the final equation above that a permuted middle structure which doesn't itself permute terminals can be obtained with an expression of the form $\underline{\underline{pxp^{-1}}}$, where $p^{-1}$ satisfies $\underline{\underline{pp^{-1}}}=I_{\mbox{\tiny{$M$}}}$.

With a more compact notation for middle structures which don't permute terminals:
\begin{equation}
(x_0, x_1, \cdots, x_n)_{\mbox{\tiny{$M$}}} \ \equiv
\ \{\Diamond x_0 \Diamond, 1\Diamond x_1 \Diamond 1, n\Diamond x_n \Diamond n\}
\end{equation}
constituent structures with multiple parallel branches can be built:
\begin{equation}
\underline{\underline{(a, b)_{\mbox{\tiny{$M$}}}(c, d)_{\mbox{\tiny{$M$}}}}} = (ac, bd)_{\mbox{\tiny{$M$}}}.
\end{equation}

\subsection{Converting Terminals to Elements of a Set}

The ordered tuple containing the empty set in each position, $(\varnothing, \varnothing, \varnothing)=\{\Diamond, \Diamond 1, \Diamond 2\}$, is a top structure which can be fused to a bottom structure, $\{\Diamond x, 1\Diamond y, 2\Diamond z\}$, to obtain a set with the corresponding elements:
\begin{equation}
\underline{\underline{(\varnothing, \varnothing, \varnothing)\{\Diamond x, 1\Diamond y, 2\Diamond z\}}} = \{x, y, z\}
\end{equation}
and, conversely, the ordered triple $(x, y, z)$ can be converted into the set $\{x, y, z\}$ by fusing it to a hypothetical bottom structure whose terminals are all empty sets:
\begin{equation}
\underline{\underline{(x, y, z)\{\Diamond , 1\Diamond , 2\Diamond \}}} = \{x, y, z\}.
\end{equation}

The set $\{\Diamond , 1\Diamond , 2\Diamond \}$ does not satisfy the definition of a bottom structure because $\Diamond$ and $1\Diamond$ are not maximal within it since they are constituents of $2\Diamond$ as well as themselves and the full set, but equation \ref{fusenotation} gives a well-defined result when $B_n=\varnothing$ for every value of $n$.

A set containing as elements multiple sets which were built in parallel using middle structures can then be obtained by performing both of these operations on a middle structure built by fusing others together. Rather than writing the full expression containing natural numbers and $\Diamond$ symbols each time, we can denote the set that results from this procedure using a notation which explicitly indicates the completion of the fusing process as well as the bracketing together of all the sets involved:
\begin{equation}
\underline{\underb{$(a, b)_{\mbox{\tiny{$M$}}}(c, d)_{\mbox{\tiny{$M$}}}$}} = \{ac, bd\}.
\end{equation}

The resulting constituent algebra is powerful enough to represent arbitrarily complicated sets in terms of constituent replacement, including the simple $\Diamond$ set:
\begin{equation}
\Diamond = \underline{\underb{$(2, 2_{{}_{_V}})_{\mbox{\tiny{$M$}}}$}}
\end{equation}
where $2=1(1)$ is Zermelo's construction of the number 2 and $2_{{}_{_V}}=\underline{\underb{$(0, 1)_{\mbox{\tiny{$M$}}}$}}$ is von Neumann's construction of the same number.

\section{Summary}

All of the objects we deal with in mathematics have constituents, while sets also have elements and cardinality. When we construct mathematical objects from sets, we should therefore encode the constituent structures of the objects in the constituent structures of the sets, not partially there and partially in the specific details of their elements and cardinality. 

If we need to refer to the cardinality or elements of the underlying set which represents an object in order to know something about that object, then we can never let go of the arbitrary choices made during the construction of the object. The properties of the object will be inseparable from the properties of the underlying set, and our representation of the object will be unnecessarily confusing and complicated. 

Whether two distinct sets constructed by two distinguished mathematicians represent the same object in the same branch of mathematics, or, by chance, objects with the same structure in different branches, is determined by the existence or non-existence of a constituent structure isomorphism between the two sets, which may have different cardinalities.

One of the principles of set theory is that a set contains only one instance of each of its elements. A set with two identical elements has one element. The concept of a constituent extends this removal of redundancy beyond a set's top layer. A set can have many copies of another set inside it, but from the point of view of functions, subsets and elements, no repetition is detectable. The constituent structure of the set is inaccessible when there's no way to express the fact that those repeated instances are all the same thing.

When we require our constructions of mathematical objects as sets to encode the object solely in the set's constituent structure, the resulting encoding is simple, natural, matches our conceptual understanding of the object, and contains a natural representation of the object's operations and relations as operations and relations on sets. 

The examples that we have covered here include natural numbers with the binary operation of addition and ordered pairs and tuples of arbitrary sets. In Part Two, we will consider more complicated but familiar mathematical objects with richer structures, specifically the integers, arithmetic expressions involving subtraction, multiplication and division, and the rational numbers.

\raggedbottom

\pagebreak

\pagebreak

\section*{Acknowledgements}

The author is deeply grateful to Venkata Rayudu Posina, whose patient explanations of category theory made the need for a category of sets modulo cardinalities visible, and also to David Marin Roma for useful discussions and suggestions which improved the readability of the article. VRP also articulated the significance for the category of constituent structures of the product and sum constructions in terms of universal properties.

\end{document}